\documentclass[12pt]{amsart}
\usepackage{amsmath,amssymb,amsbsy,amsfonts,latexsym,amsopn,amstext,cite,
                                               amsxtra,euscript,amscd,bm}
\usepackage{url}

\usepackage{mathrsfs}
\usepackage{enumitem}
\usepackage{color}
\usepackage[colorlinks,linkcolor=blue,anchorcolor=blue,citecolor=blue,backref=page]{hyperref}
\usepackage{color}
\usepackage{graphics,epsfig}
\usepackage{graphicx}
\usepackage{float}
\usepackage{epstopdf}
\hypersetup{breaklinks=true}

\usepackage[np]{numprint}
\npdecimalsign{\ensuremath{.}}

\usepackage{bibentry}

\usepackage[english]{babel}
\usepackage{mathtools}
\usepackage{todonotes}
\usepackage{url}
\usepackage[colorlinks,linkcolor=blue,anchorcolor=blue,citecolor=blue,backref=page]{hyperref}

\usepackage[norefs,nocites]{refcheck}

\theoremstyle{plain}
\newtheorem{theorem}{Theorem}
\newtheorem{lemma}[theorem]{Lemma}

\def \F{\mathbb{F}}

\def\cF{{\mathcal F}}
\def\cG{{\mathcal G}}
\def\cH{{\mathcal H}}

\def\cN{{\mathcal N}}

\def\cX{{\mathcal X}}
\def\cY{{\mathcal Y}}
\def\cZ{{\mathcal Z}}

\numberwithin{equation}{section}
\numberwithin{theorem}{section}

\def\mand{\qquad \mbox{and} \qquad}
\def\({\left(}
\def\){\right)}

\def\ep{{\mathbf{\,e}}_p}

\begin{document}

\title[Exponential sums over  subgroups and intervals]{New estimates for exponential sums over multiplicative subgroups and intervals in prime fields}
\author[D. Di Benedetto \it{et al.}]{Daniel Di Benedetto}
\address{Department of Mathematics, University of British Columbia, 1984 Mathematics Road, Vancouver, BC, V6T 1Z2, Canada}
\email{dibenedetto@math.ubc.ca}

\author[]{Moubariz Z. Garaev}
\address{Centro de Ciencias Matem{\'a}ticas, Universidad Nacional Aut{\'o}noma de M{\'e}xico, C.P. 58089, Morelia, Michoac{\'a}n, M{\'e}xico}
\email{garaev@matmor.unam.mx}

\author[]{Victor C. Garcia}
\address{Departamento de Ciencias B\'asicas, Universidad Aut\'onoma Metropo\-litana –Azcapotzalco, C.P. 02200, M\'exico D.F., M\'exico}
\email{vcgh@correo.azc.uam.mx}

\author[]{Diego Gonzalez-Sanchez}
\address{Departamento de Matem{\'a}ticas, Universidad Aut\'onoma de Madrid, and ICMAT, Madrid 28049, Spain}
\email{diego.gonzalez@icmat.es}

\author[] {Igor E. Shparlinski}
\address{Department of Pure Mathematics, University of New South Wales,
Sydney, NSW 2052, Australia}
\email{igor.shparlinski@unsw.edu.au}

\author[]{Carlos A. Trujillo}
\address{Departamento de Matem\'aticas, Universidad del Cauca
Popay\'an, Co\-lombia}
\email{trujillo@unicauca.edu.co}

%\date{}                                           % Activate to display a given date or no date

\begin{abstract}
Let $\cH$ be a multiplicative subgroup of $\mathbb{F}_p^*$ of order $H>p^{1/4}$. We show that
$$
\max_{(a,p)=1}\left|\sum_{x\in \cH} \ep(ax)\right| \le H^{1-31/2880+o(1)},
$$
where $\ep(z) = \exp(2 \pi i z/p)$,
which improves a result of Bourgain and Garaev (2009). We also obtain new estimates
for double exponential sums with product $nx$ with $x \in \cH$ and $n \in \cN$ for a short interval
$\cN$
of consecutive integers.
\end{abstract}

\maketitle
%\section{}
%\subsection{}

\section{Introduction}

 Let $p$ be a large prime number and $\F_p$ be the prime field of order $p$.
 We always assume that the elements of $\F_p$  are represented by the set $\{0,1,\ldots, p-1\}$.
 Let $\cH$ be a multiplicative subgroup of $\F_p^*$
 of order $H =|\cH|$.

 We also denote
 $$
\ep(z) = \exp(2 \pi i z/p).
$$

The problem of obtaining nontrivial upper bounds for the exponential sum
\begin{equation}
\label{eq:SaH}
S_a(\cH) = \sum_{x\in \cH}\ep(a x)
\end{equation}
is a classical problem with a variety of results and applications in number theory.   The classical result of  Gauss implies  that if $H = (p-1)/2$, then $|S_a(\cH)|=p^{1/2}$. From the work of Hardy and Littlewood on the Waring problem  it is known
that $|S|<p^{1/2}$, which is non-trivial when $H > p^{1/2}$. The problem of obtaining nontrivial bounds for $H < p^{1/2}$ has been a subject of much research, see~\cite{HBK, Kon, Shp, BGK}.
Using the sum-product estimate and other tools from additive combinatorics Bourgain, Konyagin and Glibichuk~\cite{BGK}  proved that if $H>p^{\varepsilon}$,
then
$$
\left| S_a(\cH) \right| \le H  p^{-\delta},
$$
where $\delta>0$ depends only on $\varepsilon$.  Prior to their work, this
estimate had been only known under the assumption
$H >p^{1/4+\varepsilon}$ due to Konyagin. In the limiting case $H \sim p^{1/4}$ Bourgain and Garaev~\cite{BG2009} obtained the following  explicit bound:
\begin{equation}
\label{eqn:BG2009}
\max_{(a,p)=1}\left| S_a(\cH) \right| \le H^{1-175/9437184+o(1)}.
\end{equation}
The argument of~\cite{BG2009} is based on explicit trilinear exponential sum estimates obtained in the same paper, and also on a bound of Konyagin~\cite{Kon} on the number $T_m(\cH)$
of solutions of the congruence
\begin{equation}
\label{eqn:defTm}
  h_1 + \cdots + h_m \equiv h_{m+1} + \cdots + h_{2m} \pmod p, \quad  h_1,\ldots, h_{2m} \in \cH.
\end{equation}

Since the work~\cite{BG2009} there have been significant developments on sum-product problems which have lead to new trilinear sum estimates of Macourt~\cite{Mac} and Petridis and Shparlinski~\cite{PS2019}.
Moreover, new bounds for $T_m(\cH)$ have been obtained by Murphy, Rudnev, Shkredov and
Shteinikov~\cite[Theorem~3 and Corollary~7]{MRSS2017} for the cases $m=2$ and $m=3$, and by Shkredov~\cite{S2018} for the case $m\ge 4$. In the present paper, combining these estimates with the argument from~\cite{BG2009},
we improve the estimate~\eqref{eqn:BG2009}, replacing
$$
175/9437184 = 1.854 \ldots \times 10^{-5}
\quad \text{with} \quad 31/2880 = 1.076 \ldots \times 10^{-2}.
$$

Next, we consider the double sum involving intervals and subgroups. Let
$$
\cN = \{L+1,\ldots, L+N\}
$$
be an interval of consecutive $N$ integers with $|\cN|=N\le p$. For $\gcd(a,p)=1$, we denote
\begin{equation}
\label{eq:SaNH}
S_a(\cN, \cH)=\sum_{n\in \cN}\left|\sum_{x\in \cH}\ep(anx)\right|.
\end{equation}

This sum is a special case of a more general family of exponential sums considered by 
Bourgain~\cite{Bour}, and more recently
by Garaev~\cite{G2019} and Shparlinski and Yau~\cite{ShpYau}.  In the present paper we obtain new estimates for
 $S_{a}(\cN,\cH)$.

\section{Notation}

 In what follows, we use the notation $A\lesssim B$ to mean that
$|A|<B p^{o(1)}$, or equivalently,  for any $\varepsilon>0$ there
is a constant $c(\varepsilon)$, which depends only on $\varepsilon$, such that
$|A| \le c(\varepsilon) Bp^{\varepsilon}$.
We also recall that  the standard notations $A=O(B)$, $A\ll B$ and $B \gg A$ are each equivalent to the
statement that the inequality $|A|\le c\,B$ holds with a
constant $c>0$ which is absolute throughout this paper.

 \section{Our results}
We start with a bound on the sums $S_a(\cH)$ over small subgroups given by~\eqref{eq:SaH}.

\begin{theorem}
\label{thm:SingleSum}
Let $\cH$ be a multiplicative subgroup of $\F_p^*$ of order $H$ with $p^{1/2} > H>p^{1/4}$. Then
    \[\max_{(a,p)=1}\left|S_a(\cH) \right| \lesssim  H^{2689/2880} p^{1/72}.\] 
\end{theorem}

In particular, when $H>p^{1/4}$,  Theorem~\ref{thm:SingleSum} gives
\[\max_{(a,p)=1}\left|S_a(\cH) \right| \lesssim H^{1-31/2880}.\]

Next we consider the sums $S_a(\cN, \cH)$ over small subgroups given by~\eqref{eq:SaNH}.
It is convenient to define
$$
\Gamma(\cN, \cH) = 1 + \frac{H}{N}+\frac{NH}{p}+\frac{H^{3/4}}{p^{1/4}}.
$$
It is useful to observe that for  $H \lesssim N  \lesssim p^{1/3}$ we have 
$\Gamma(\cN, \cH)   \lesssim 1$. 

\begin{theorem}
\label{thm:S(a) 2 3} 
Let $\cH$ be a multiplicative subgroup of $\F_p^*$ of order $H<p^{1/2}$ and 
 and   let $\cN$   be an interval of consecutive $N$   integers. Then
\[
S_a(\cN, \cH) \lesssim NH\times\min\{\Delta_1^{1/4}, \Delta_2^{1/6}\},
\]
where
$$
\Delta_1 = \frac{p}{NH^{2+ 11/20}} \Gamma(\cN, \cH)
\mand
 \Delta_2 = \frac{p}{NH^3}\Gamma(\cN, \cH) .
$$
\end{theorem}

In particular, if $N, H = p^{1/3+o(1)}$, then we get the bound
$$
S_a(\cN, \cH) \lesssim N^{2- 1/6}.
$$
This improves the result of~\cite[Theorem~1]{G2019} in the case of subgroups.

Note that Theorem~\ref{thm:S(a) 2 3} 
is based on results of~\cite[Theorem~3 and Corollary~7]{MRSS2017} and is trivial when $NH^3 < p$.   
By using a result of~\cite{S2018} in the proof instead, one can improve this range. 
However, proceeding as in the proof of Theorem~\ref{thm:SingleSum}, we also obtain the following stronger result.

\begin{theorem}\label{thm:S(N,H)short}
  Let $\cH$ be a multiplicative subgroup of $\F_p^*$ of order $H<p^{1/2}$ and 
   let $\cN$   be an interval of consecutive $N$ integers  which avoids $0$ modulo $p$.   Then
$$
     S_a(\cN, \cH) \lesssim NH  \Delta^{1/24} 
$$
    where
  \[
     \Delta=\frac{p}{N H^{3+31/40}} \Gamma(\cN,\cH).
  \]  
\end{theorem}

   We note that in the case when $N=H$,  Theorem~\ref{thm:S(N,H)short}   together with Theorem~\ref{thm:S(a) 2 3}, gives an 
   improvement to the bound resulting from a direct application of Theorem~\ref{thm:SingleSum} 
   in the full range for which Theorem~\ref{thm:SingleSum} is nontrivial. In particular, 
   if $N,H = p^{1/4+o(1)}$ we obtain 
\[
   S_a(\cN, \cH) \lesssim N^{2- 31/960},
   \]   
and observe that  
\[
31/960 =3.229\ldots \times 10^{-2} > 31/2880 = 1.076 \ldots \times 10^{-2}.
\]

Furthermore, it is easy to see that for $H\sim p^{1/4}$  Theorem~\ref{thm:S(N,H)short} improves 
the bound $NH^{1-31/2880}$ which follows directly from  Theorem~\ref{thm:SingleSum}
provided that $N \ge p^{89/480 + \varepsilon}$ for some fixed $\varepsilon>0$.

%--------------

\section{Tools}

We need the following trilinear exponential sum bound, which is due to Petridis and Shparlinski~\cite[Theorem~1.1]{PS2019}.

\begin{lemma}
\label{tri}
For any sets $\cX,\cY,\cZ\subseteq\mathbb{F}_p^*$ and any complex numbers $\alpha_x,\beta_y,\gamma_z$ with $|\alpha_x|\leq1,|\beta_y|\leq1,|\gamma_z|\leq1$ we have
    \[\sum_{x\in \cX}\sum_{y\in \cY}\sum_{z\in \cZ} \alpha_x \beta_y \gamma_z \ep(axyz)\ll p^{1/4}|\cX|^{3/4}|\cY|^{3/4}|\cZ|^{7/8}.\]
\end{lemma}

We recall that $\cH$ is a multiplicative subgroup of $\F_p^*$ with $H=|\cH|$ elements, and  $T_m(\cH)$  the number of solutions of the congruence~\eqref{eqn:defTm}.
The following two results are due to Murphy,  Rudnev, Shkredov and  Shteinikov~\cite[Theorem~3 and Corollary~7]{MRSS2017}.

\begin{lemma}
\label{lem: T2 energy} Let  $H < \sqrt{p}$. Then
\[
 T_2(\cH)\ll H^{49/20}\log^{1/5}H.
\]
\end{lemma}

\begin{lemma}\label{lem:tri-sum-energy} Let  $H < \sqrt{p}$. Then
\[
 T_3(\cH)\ll H^4\log H.
\]
\end{lemma}

Let now $J(\cN,\cH)$ be the number of solutions of the equation
\[
 n_1 h_1 \equiv n_2 h_2 \pmod p, \quad   n_1,n_2 \in \cN, \ h_1,h_2 \in \cH.
\]
 We also need the  result from~\cite[Corollary 1]{G2019} which is based on the work of Cilleruelo and Garaev~\cite{CilGar}.

\begin{lemma}\label{lem:Cill-Gar-bound}
The following bound holds:
\[
  J(\cN,\cH) \lesssim H^2 + NH + \frac{N^2H^2}{p}+\frac{NH^{7/4}}{p^{1/4}}.
\]
\end{lemma}
We remark that when $\cN$ starts from the origin (that is, if $L=0$ in the definition of~$\cN$), 
the first term $H^2$ on the right hand side can be removed.
It is also to be mentioned that a result similar to Lemma~\ref{lem:Cill-Gar-bound} has been obtained 
by Banks and Shparlinski~\cite{BSh} for arbitrary sets instead of $\cH$.

\section{Proof of Theorem~\ref{thm:SingleSum}}

Fix an element $a\in\mathbb{F}_p^*$ and let
    \[\Delta = \frac{1}{H} \left| \sum_{x\in \cH} \ep(ax) \right|.\]
Hence we  now have to prove an upper bound on $\Delta$. 
We can assume that 
\begin{equation}
\label{eq: large D}
\Delta>H^{-37/960}
\end{equation}
which the largest saving we have in the bound of Theorem~\ref{thm:SingleSum}
when $H$ is approaching $p^{1/2}$.

Clearly, since $H$ is a multiplicative subgroup, for any $y\in \cH$ we have
\begin{align*}
S_a(\cH)^3 &= \sum_{x_1,x_2,x_3\in \cH}\ep\(a(x_1+x_2+x_3)\)\\
& = \sum_{x_1,x_2,x_3\in \cH} \ep(a(x_1+x_2+x_3)y).
\end{align*} 
Hence
$$
H S_a(\cH)^3=  \sum_{y\in \cH} \sum_{x_1,x_2,x_3\in \cH}  \ep\(a(x_1+x_2+x_3)y\).
$$
Recalling that $| S_a(\cH)| = H\Delta$ and changing the order of summation, we obtain 
\begin{equation}
\begin{split}
\label{eq: H3|H|}
H^{4}\Delta^3  
& =  \left |  \sum_{x_1,x_2,x_3\in \cH}   \sum_{y\in \cH} \ep\(a(x_1+x_2+x_3)y\)\right|\\
& \le   \sum_{x_1,x_2,x_3\in \cH}    \left | \sum_{y\in \cH} \ep\(a(x_1+x_2+x_3)y\)\right|.
\end{split}
\end{equation}

Let $\cF^*$ be the set of triples $(x_1,x_2,x_3)\in \cH^3$ with
$$
 \left | \sum_{y\in \cH} \ep\(a(x_1+x_2+x_3)y\)\right| \ge \frac{1}{2} H \Delta^3 .
$$
Discarding the contribution to the right hand side of~\eqref{eq: H3|H|}
(not exceeding $0.5 H^4 \Delta^3$)
from the triples  $(x_1,x_2,x_3)\notin \cF^*$, we have
\begin{equation}
\label{eq: F|H|}
\sum_{(x_1,x_2,x_3)\in \cF^*}    \left | \sum_{y\in \cH} \ep\(a(x_1+x_2+x_3)y\)\right| \ge  \frac{1}{2} H^4 \Delta^3.
\end{equation}

Let $\cF_i$ be the set of triples $(x_1,x_2,x_3)\in \cF^*$ with
$$
H2^{-i-1}  <  \left | \sum_{y\in \cH} \ep\(a(x_1+x_2+x_3)y\)\right| \le H 2^{-i}, \qquad i =0,1, \ldots.
$$

Clearly, the sets $\cF_i$ are non-empty for $O(\log H)$ values of $i$. 
Thus, we see from~\eqref{eq: F|H|} that  there exists $i_0 = O(\log H)$ with
\begin{align*}
H 2^{-i_0} | \cF_{i_0} |& \ge  \sum_{(x_1,x_2,x_3)\in \cF_{i_0} }   \left | \sum_{y\in \cH} \ep\(a(x_1+x_2+x_3)y\)\right|\\
& \gg \frac{1}{\log p}   \sum_{x_1,x_2,x_3\in \cH}    \left | \sum_{y\in \cH} \ep\(a(x_1+x_2+x_3)y\)\right|\\
& \ge \frac{1}{\log p} H^{4}\Delta^3.
\end{align*}
Setting $\cG_1 =   \cF_{i_0}$, 
we see that  we can find a set $\cG_1\subseteq \cH^3$ and a number $\Delta_1$,
with
\begin{equation}
\label{eq: Set G1}
|\cG_1|\gtrsim \frac{H^3\Delta^3}{\Delta_1}  \mand \Delta^3 \lesssim \Delta_1< 1,
\end{equation}
such that for any $(x_1,x_2,x_3)\in \cG_1$ we have
\begin{equation}
\label{eq: |H|}
\left|\sum_{y\in \cH }\ep(a(x_1+x_2+x_3)y)\right| \gg H\Delta_1.
\end{equation}

Let $\cX=\{x_1+x_2+x_3:~(x_1,x_2,x_3)\in \cG_1\}\setminus\{0\}$.
To estimate the cardinality $|\cX|$ we denote by $J(x)$ the number of representations
of $x \in \cX\cup\{0\}$ as $x = x_1+x_2+x_3$, $(x_1,x_2,x_3)\in \cG_1$.
Clearly
$$
\sum_{x \in \cX\cup\{0\}} J(x) = |\cG_1| \mand \sum_{x \in \cX\cup\{0\}} J(x)^2 \le T_3(\cH).
$$
Applying the Cauchy--Schwarz inequality, we see that by Lemma~\ref{lem:tri-sum-energy} and the bound~\eqref{eq: Set G1} we have
$$
|\cX\cup\{0\}|\ge \frac{|\cG_1|^2}{ T_3(\cH)} \gtrsim H^2\Delta^6/\Delta_1^2.
$$
Recalling our assumption~\eqref{eq: large D} we see that if $p$ is large enough then 
\begin{equation}
\label{eq: large X}
|\cX| \gtrsim H^2\Delta^6/\Delta_1^2, 
\end{equation}
and summing~\eqref{eq: |H|} over $x\in \cX$ we obtain
\[\sum_{x\in \cX} \left|\sum_{y\in \cH }\ep(axy)\right| \geq H|\cX|\Delta_1.\]
Applying the H\"older inequality gives
\begin{equation}
\begin{split} 
\label{eq: Holder}
\sum_{x\in \cX} \left|\sum_{y_1,y_2,y_3\in \cH}\ep\(ax(y_1+y_2+y_3)\)\right| & \geq
\frac{\(H|\cX|\Delta_1\)^3}{|\cX|^2}\\
& = H^3|\cX|\Delta_1^3.
\end{split}
\end{equation}
Since $\cH$ is a multiplicative subgroup, for any $z\in \cH$ we have
$$
\sum_{y_1,y_2,y_3\in \cH}\ep\(ax(y_1+y_2+y_3)\)
 =\sum_{y_1,y_2,y_3\in \cH}\ep\(ax(y_1+y_2+y_3)z\).
 $$
 Therefore, we derive from~\eqref{eq: Holder} that
 \begin{align*}
H^4|\cX|\Delta_1^3  &
=  \sum_{x\in \cX} \left|\sum_{y_1,y_2,y_3\in \cH}  \sum_{z\in \cH} \ep\(ax(y_1+y_2+y_3)z\)\right|\\
& \le  \sum_{x\in \cX} \sum_{y_1,y_2,y_3\in \cH}  \left| \sum_{z\in \cH} \ep\(ax(y_1+y_2+y_3)z\)\right|.
\end{align*}

We now repeat the previous dyadic argument with respect to triples $(y_1,y_2,y_3)\in \cH^3$,
and see that we can again find a set $\cG_2\subseteq \cH^3$ and a number $\Delta_2$,
with
\begin{equation}
\label{eq: Set G2}
|\cG_2|\gtrsim \frac{H^3\Delta_1^3}{\Delta_2} \mand \Delta_1^3\lesssim \Delta_2\le 1
\end{equation}
such that for any $(y_1,y_2,y_3)\in \cG_2$ we have
\begin{equation}
\label{eq: X|H|}
\sum_{x\in \cX}\left|\sum_{z\in \cH }\ep(ax(y_1+y_2+y_3)z)\right| \geq H|\cX| \Delta_2.
\end{equation}

Let $\cY=\{y_1+y_2+y_3 :~(y_1,y_2,y_3)\in \cG_2\}\setminus\{0\}$.
Similarly to~\eqref{eq: large X}, by Lemma~\ref{lem:tri-sum-energy} and
the bound~\eqref{eq: Set G2}, we have
\begin{equation}
\label{eq: large Y}
|\cY|\gtrsim H^2\Delta_1^6/\Delta_2^2,
\end{equation}
 and summing~\eqref{eq: X|H|} over $y\in \cY$ we obtain
   \[\sum_{y\in \cY}\sum_{x\in \cX}\left|\sum_{z\in \cH }\ep(axyz)\right| \geq H|\cX||\cY|\Delta_2.\]
Applying the Cauchy--Schwarz inequality, we get
    \[ \sum_{z_1,z_2\in \cH}\left|\sum_{x\in \cX}\sum_{y\in \cY} \ep(axy(z_1-z_2))\right| \geq H^2|\cX||\cY|\Delta_2^2.\]

Again from the previous dyadic argument however with respect to pairs $(z_1,z_2) \in  \cH^2$,
we see that 
we can find a set $\cG_3\subseteq \cH^2$ and a number $ \Delta_3$  with
\begin{equation}
\label{eq: Set G3}
|\cG_3| \gtrsim \frac{H^2\Delta_2^2}{\Delta_3} \mand \Delta_2^2\lesssim \Delta_3\le 1
\end{equation}
such that for any $(z_1,z_2)\in \cG_3$ we have
\begin{equation}
\label{eq: XY|H|}
\left|\sum_{x\in \cX}\sum_{y\in \cY} \ep(axy(z_1-z_2))\right| \geq  |\cX||\cY| \Delta_3.
\end{equation}

Let $\cZ=\{z_1-z_2 :~(z_1,z_2)\in \cG_3\}\setminus\{0\}$. Similarly to~\eqref{eq: large X} and~\eqref{eq: large Y}, by Lemma~\ref{lem: T2 energy} and
the bound~\eqref{eq: Set G3}, we have
\begin{equation}
\label{eq: large Z}
|\cZ| \gtrsim H^{31/20}\Delta_2^4/\Delta_3^2,
\end{equation}
and summing~\eqref{eq: XY|H|}  over $z\in \cZ$ we obtain
    \[\sum_{z\in \cZ}\left|\sum_{x\in \cX}\sum_{y\in \cY}  \ep(axyz)\right| \geq |\cX||\cY||\cZ| \Delta_3.\]
Comparing this with Lemma~\ref{tri} implies that
\[
  |\cX||\cY||\cZ|^{1/2}\Delta_3^4 \ll p.
\]
Finally, by applying the estimates~\eqref{eq: large X}, \eqref{eq: large Y}
and~\eqref{eq: large Z}  we conclude that
\begin{equation}
\label{eq: penult}
p \gtrsim \frac{H^2\Delta^6}{\Delta_1^2} \cdot  \frac{ H^2\Delta_1^6}{\Delta_2^2}  \cdot  \frac{H^{31/40}\Delta_2^2}{\Delta_3}
\Delta_3^4  = H^{191/40}  \Delta^6  \Delta_1^4  \Delta_3^3.
\end{equation}
Since by~\eqref{eq: Set G1}, \eqref{eq: Set G2} and~\eqref{eq: Set G1} we have
\begin{align*}
 \Delta^6  \Delta_1^4  \Delta_3^3  & \gtrsim
 \Delta^6  (\Delta^3)^4  \(\Delta_2^2\)^3  =  \Delta^{18}  \Delta_2^6\\
 & \gtrsim  \Delta^{18}  \(\Delta_1^3\)^6 =   \Delta^{18}  \Delta_1^{18} \\
 & \gtrsim  \Delta^{18}  \(\Delta^3\)^{18}  =  \Delta^{72}
\end{align*}
together with~\eqref{eq: penult}, we obtain
\[
  \Delta \lesssim p^{1/72}H^{-191/2880}
\]
and conclude the proof.

\section{Proof of Theorem~\ref{thm:S(a) 2 3}}

We use the abbreviation $S = S_a(\cN, \cH)$. Since $H$ is a subgroup of $\F_p^*$, for every $u\in \cH$ we have that
\[
   S = \sum_{n\in \cN} \left|\sum_{h\in \cH}\ep(anhu)\right|.
\]
Summing over $u\in \cH$ we get
\[
  S= \frac{1}{H} \sum_{n\in \cN} \sum_{h\in \cH}\left|\sum_{u\in \cH}\ep(anhu) \right| .
\]
Let $m \ge 2$ be an integer constant. Applying the H\"older inequality, we obtain
$$
  S^m\le  \frac{N^{m-1}}{H}\sum_{n\in \cN}\sum_{h\in \cH}\left|\sum_{u\in \cH}\ep(anhu)\right|^m.
$$

Denote by $R(\lambda)$ be the umber of solutions of the congruence
$$
nh\equiv \lambda\pmod p,\quad n\in \cN, \,h\in\cH.
$$
It follows that
$$
  S^m \le  \frac{N^{m-1}}{H}\sum_{\lambda=0}^{p-1}R(\lambda)\left|\sum_{u\in \cH}\ep(a\lambda u)\right|^m.
$$
Squaring out and applying the Cauchy--Schwarz inequality, we obtain that
$$
S^{2m} \le \frac{N^{2m-2}}{H^2} \sum_{\lambda=0}^{p-1}R^2(\lambda) \sum_{\mu=0}^{p-1}\left|\sum_{u\in \cH}\ep(a\mu u)\right|^{2m}.
$$
Clearly,
$$
\sum_{\lambda=0}^{p-1}R^2(\lambda) = J(\cN, \cH) \mand \sum_{\mu=0}^{p-1}\left|\sum_{u\in \cH}\ep(a\mu u)\right|^{2m} =p\,T_{m}(\cH).
$$
Hence, applying Lemma~\ref{lem:Cill-Gar-bound}, we get that
$$
S^{2m}\lesssim \frac{N^{2m-2}}{H^2}\(NH +H^2 + \frac{N^2H^2}{p}+\frac{NH^{7/4}}{p^{1/4}}\) p\, T_{m}(\cH).
$$
Therefore,
$$
S\lesssim NH \Delta^{1/2m},
$$
where
\[
\Delta = \frac{p\,T_m(\cH)}{NH^{2m+1}} \Gamma(\cN, \cH).
\]
Consequently taking $m=2,3$ and applying Lemmas~\ref{lem: T2 energy} and~\ref{lem:tri-sum-energy}, we obtain the desired result.

\section{Proof of Theorem~\ref{thm:S(N,H)short}}

  As mentioned, we will follow the strategy of the proof of Theorem~\ref{thm:SingleSum}.   Set
\[
   \Delta=\frac{1}{NH}  \sum_{n \in \cN}\left| \sum_{h \in \cH} \ep(anh) \right|.
\]   
  Besides, we suppose that $\Delta > H^{-37/320}$ (otherwise the result is trivial).
  Since $\cH$ is a subgroup, for any $u,v \in \cH$ we have
\[
   NH \Delta = \sum_{n \in \cN}\left| \sum_{h \in \cH} \ep(anhuv) \right|,
\]  
  and the sum does not depend on $u$ or $v$. Then
\begin{align*}
  H^3 N \Delta & =  \sum_{u \in \cH} \sum_{v \in \cH} \sum_{n \in \cN}\left| \sum_{h \in \cH} \ep(anhuv) \right|\\
    & \le \sum_{n \in \cN} \sum_{h \in \cH} \left| \sum_{u \in \cH} \sum_{v \in \cH} \ep(anhuv) \right|.
\end{align*}
  There is a subset $\cG_1 \subseteq\cN \times \cH$ and $\Delta_1 >0$, with
\[  
  |\cG_1|\gtrsim N H \frac{\Delta}{\Delta_1} \mand  \Delta \lesssim \Delta_1 <1, 
\]
  such that for any $(n,h) \in \cG_1$ we have
\begin{equation}\label{ineq:NH}
    \left| \sum_{u \in \cH} \sum_{v \in \cH} \ep(anhuv) \right| \ge  H^2 \Delta_1.
\end{equation}
  Define $\cX=\{nh:~ (n,h)\in \cG_1\}$ and note that by the hypothesis on $\cN$, 
  it cannot contain $0$. 
  By Lemma~\ref{lem:Cill-Gar-bound},
\begin{equation}\label{ineq:X}
  |\cX| \ge \frac{|\cG_1|^2}{J(\cN,\cH)}\gtrsim \frac{NH}{\Gamma(\cN,\cH)} 
  \frac{\Delta^2}{\Delta_1^2}.
\end{equation}

   Then, for each $x \in \cX$ satisfying~\eqref{ineq:NH}, we get
\begin{align*}
   H^6 \Delta_1^3 & \le  \left|\sum_{u \in \cH} \sum_{v \in \cH} \ep(axuv) \right|^3 \le 
       H^2 \sum_{u \in \cH} \left|\sum_{v \in \cH}  \ep(axuv) \right|^3 \\
          & \le  H^2 \sum_{u \in \cH} \theta^3 \sum_{v_1, v_2, v_3 \in \cH}  \ep(axu(v_1 + v_2 + v_3))  \\
         & \le  H^2  \sum_{v_1, v_2, v_3 \in \cH} \left|\sum_{u \in \cH}  \ep(axu(v_1 + v_2 + v_3)) \right|,
\end{align*}
   where $|\theta|=1$, and does not depend on $u$. Hence, 
   we find a subset $\cG_2 \subseteq\cH^3$ and $\Delta_2 >0$, with
\[  
  |\cG_2|\gtrsim H^3 \frac{\Delta_1^3}{\Delta_2} \mand  \Delta_1^3 \lesssim \Delta_2 <1, 
\]
  such that for any $(v_1,v_2,v_3) \in \cG_2$ we have
\begin{equation}\label{ineg:V's}
    \left|\sum_{u \in \cH}  \ep(axu(v_1 + v_2 + v_3)) \right| \ge  H \Delta_2.
\end{equation}  
  Let $\cY = \{v_1+v_2+v_3 : (v_1,v_2,v_3)\in \cG_2\}\setminus \{0\},$ and 
  applying Lemma~\ref{lem:tri-sum-energy} we have 
\[
  |\cY \cup \{0\}| \ge \frac{|\cG_2|^2}{T_3(\cH)}\gtrsim H^2 \frac{\Delta_1^6}{\Delta_2^2}.
\]
  Recalling the assumption $\Delta > H^{-37/320},$ we note that if $p$ is large enough then
\begin{equation}\label{ineq:Y}
  |\cY| \gtrsim H^2 \frac{\Delta_1^6}{\Delta_2^2}.
\end{equation}
  Given $x\in \cX$ and $y \in \cY,$ then squaring out~\eqref{ineg:V's} we get
\[
  H^2 \Delta_2^2 \le \left|\sum_{u \in \cH} \ep(a x y u) \right|^2 
        = \sum_{u_1, u_2 \in \cH}\ep(axy (u_1-u_2)).
\]   
  Summing over $x \in \cX$ and $y\in \cY$ it follows that
\begin{align*}
  H^2 |\cX| |\cY| \Delta_2^2 & \le \sum_{x \in \cX}\sum_{y \in \cY}\sum_{u_1, u_2 \in \cH}\ep(axy (u_1-u_2)) \\
      & \le \sum_{u_1, u_2 \in \cH} \left| \sum_{x \in \cX}\sum_{y \in \cY} \ep(axy(u_1-u_2)) \right|.
\end{align*}
  Now again, we find a subset $\cG_3 \subseteq \cH^2$ and $\Delta_3 >0$, with
\[  
  |\cG_3|\gtrsim H^2 \frac{\Delta_2^2}{\Delta_3} \mand  \Delta_2^2 \lesssim \Delta_3 <1, 
\]
  such that for any $(u_1,u_2) \in \cG_3$ we have
\begin{equation}\label{ineq:U's}
    \left|\sum_{x \in \cX} \sum_{y \in \cY}  \ep(axy(u_1 - u_2)) \right| \ge  |\cX||\cY| \Delta_3.
\end{equation}  
  Let $\cZ=\{u_1-u_2\;:\; (u_1,u_2)\in \cG_3\} \setminus \{0\}$, and applying 
  Lemma~\ref{lem: T2 energy} we have 
\[
   |\cZ \cup \{0\}| \ge \frac{|\cG_3|^2}{T_2(\cH)}\gtrsim {H^{31/20}} \frac{\Delta_2^4}{\Delta_3^2}.
\]
By an argument analogous to previous ones, we claim that if $p$ is large enough then
\begin{equation}\label{ineq:Z}
  |\cZ| \gtrsim {H^{31/20}} \frac{\Delta_2^4}{\Delta_3^2}.
\end{equation}  

  Therefore, summing~\eqref{ineq:U's} over $z \in \cZ$,  we obtain
\[
   |\cX||\cY||\cZ| \Delta_3 \le \sum_{z \in \cZ} \left|\sum_{x \in \cX} \sum_{y \in \cY}  \ep(axyz) \right|.
\]
  Comparing this with Lemma~\ref{tri}, we get
\[
   |\cX|^2|\cY|^2|\cZ| \,\Delta_3^8 \ll p^2.
\]  
  Now, by the estimates~\eqref{ineq:X}, \eqref{ineq:Y}, \eqref{ineq:Z} and the
  relations $\Delta_1 \gtrsim \Delta$, 
  $\Delta_2 \gtrsim \Delta_1^3$, $\Delta_3 \gtrsim \Delta_2^2$, it follows
  that
$$
  p^2  \gg \Delta_3^8 |\cX|^2 |\cY|^2 |\cZ| \gtrsim \frac{N^2 H^{6+31/20}}{\Gamma^2(\cN,\cH)}\Delta^4\Delta_1^8 \Delta_3^6 
     \gtrsim \frac{N^2 H^{6+31/20}}{\Gamma^2(\cN,\cH)}\Delta^{48},
$$
which concludes the proof. 

\section*{Acknowledgments}

D.DB., M.G., V.G., D.G-S. and C.T. are grateful to the organizers of the CMO-BIRS workshop ``Number Theory in the Americas", Oaxaca, 11-17 August 2019, for the excellent working conditions.
D.DB. was partially supported by graduate fellowships from the University of British Columbia. D.G-S. was partially supported by ``la Caixa" Foundation (ID 100010434), under agreement LCF/BQ/SO16/52270027.  V.G. and D.G-S. were supported by  A. V. Humboldt Foundation (via H. A. Helfgott's chair). I.S. was  supported   by ARC Grant~DP170100786.

\end{document}